\documentclass[10pt,twocolumn]{IEEEtran}
\usepackage{tikz}
\usepackage{cite}
\usepackage{graphicx}
\usepackage{amsmath}
\usepackage{amsfonts}
\usepackage{amssymb}
\usepackage{mathrsfs}
\usepackage{fix2col}
\usepackage{latexsym}
\usepackage[dvips]{epsfig}
\usepackage[font=footnotesize]{subfig}
\usepackage{hyperref}
\usepackage{graphics}
\usepackage{epstopdf}
\usepackage{enumerate}
\usepackage{mathtools, cuted}
\usepackage{lipsum, color}
\usepackage{epstopdf}
\usepackage{amsthm}
\usepackage{arydshln}
\usepackage{multirow}
\usepackage[]{algorithm2e}
\newtheorem{theorem}{Theorem}
\newtheorem{definition}{Definition}
\newtheorem{lemma}{Lemma}

\newtheorem{remark}{Remark}

\newtheorem{example}{Example}

\usepackage{color}

\begin {document}

\title{{\Large\bf A Bound Strengthening Method for Optimal Transmission Switching in Power Systems}
	\thanks{Email: \{fattahi, lavaei, atamturk\}@berkeley.edu}}
\author{
	\authorblockN{Salar Fattahi$^*$, Javad Lavaei$^{*,+}$, and Alper Atamt\"urk$^*$\\ \vspace{3mm}}
	\authorblockA{$^*$ Industrial Engineering and Operations Research,  University of California, Berkeley\\
	$^+$ Tsinghua-Berkeley Shenzhen Institute, University of California, Berkeley.}
	\thanks{This work was supported by the ONR YIP Award, DARPA YFA Award, AFOSR YIP Award, NSF CAREER
		Award, and an ARL Grant. A. Atamt\"urk was supported, in part, by grant FA9550-10-1-0168 from the Office of the Assistant Secretary of Defense for Research and Engineering.
	
Parts of this paper have appeared in the conference paper [22].}}
\maketitle

\begin{abstract}
	This paper studies the optimal transmission switching (OTS) problem for power systems, where certain lines are fixed (uncontrollable) and the remaining ones are controllable via on/off switches. The goal is to identify  a topology of the power grid that minimizes the cost of the system operation while satisfying the physical and operational constraints. Most of the existing methods for the problem are based on first converting the OTS into a mixed-integer linear program (MILP) or mixed-integer quadratic program (MIQP), and then iteratively solving a series of its convex relaxations. The performance of these methods depends heavily on the strength of the MILP or MIQP formulations. In this paper, it is shown that finding the strongest variable upper and lower bounds to be used in an MILP or MIQP formulation of the OTS based on the big-$M$ or McCormick inequalities is NP-hard. Furthermore, it is proven that unless $P=NP$, there is no constant-factor approximation algorithm for constructing these variable bounds. Despite the inherent difficulty of obtaining the strongest bounds in general, a simple bound strengthening method is presented to strengthen the convex relaxation of the problem when there exists a connected spanning subnetwork of the system with fixed lines. The proposed method can be treated as a preprocessing step that is independent of the solver to be later used for numerical calculations and can be carried out offline before initiating the solver. A remarkable speedup in the runtime of the mixed-integer solvers is obtained using the proposed bound strengthening method for medium- and large-scale real world systems.
\end{abstract}

\vspace{2mm}
\section{Introduction}
In power systems, transmission lines have traditionally been considered  uncontrollable infrastructure devices,  except in the case of an outage or maintenance. However, due to the pressing needs to boost the sustainability, reliability and efficiency, power system directors call on  leveraging the flexibility in the topology of the grid and co-optimizing the production and topology to improve the dispatch. In the last few years, Federal Energy Regulatory Commission (FERC) has held an annual conference on ``Increasing Market and Planning Efficiency through Improved Software"~\cite{FERC} to encourage research on the development of efficient software for enhancing the efficiency of the power systems via optimizing the flexible assets (e.g., transmission switches) in the system. Furthermore, The Energy Policy Act of 2005 explicitly addresses the ``difficulties of siting major new transmission facilities" and calls for the utilization of better transmission technologies~\cite{FERC_Act}.   

%Removing a few lines from the power network could increase the power transfer. This would, in theory, stand in  contrast with a conventionally accepted fact in the area of network flows; removing edges from the network would deteriorate the performance or increase the cost of the optimal flows. This seemingly inconsistent phenomenon stems from the underlying physical laws, which delineates a fundamental difference between power systems and classical network flow models. By taking these laws into account and harnessing the flexibilities in the topology of the grid, one can indeed improve the performance of  power systems \cite{Hedman2011}.

Unlike in the classical network flows, removing a line from a power network may improve the efficiency of the network due to physical laws. This phenomenon has been observed and harnessed to improve the power system performance by many authors.
%The transmission switching problem has long been studied as a corrective action for the line overloading and voltage violations~\cite{Bacher86, Rolim99, Shao2006, PJM2008, Granelli2006}. As another line of work, researchers have focused on minimizing the system loss (as opposed to the cost of operations)~\cite{Bacher88, Zaoui07}.
The notion of {\it optimally switching the lines of a transmission network} was introduced by O'Neill {\it et al.} \cite{ONeill2005}. 
Later on, it has been shown in a series of papers that the incorporation of controllable transmission switches in a grid could relieve network congestions \cite{Shao2006}, serve as a corrective action for voltage violation~\cite{Bacher86, Rolim99, Granelli2006}, reduce system loss~\cite{Bacher88, Zaoui07} and operational costs \cite{ Hedman2010}, improve the reliability of the system \cite{Hedman2009, Khanabadi2013} and enhance the economic efficiency of power markets \cite{Hedman2011_2}. 
We refer the reader to Hedman {\it et al.}~\cite{Hedman2011} for a survey on the benefits of transmission switching in power systems.
However, the identification of an optimal topology, namely  optimal transmission switching (OTS) problem, is a non-convex combinatorial optimization problem that is proved to be NP-hard \cite{Lehmann14}. Therefore,  brute-force search algorithms for finding an optimal topology are often inefficient. Most of the existing methods are based on heuristics and iterative relaxations of the problem. These methods include, but are not restricted to, Benders decomposition \cite{Hedman2010, Khanabadi2013}, branch-and-bound and cutting-plane methods \cite{Fisher2008,Burak2016}, genetic algorithms \cite{Granelli2006}, and line ranking \cite{Barrows12, Fuller2012}. Recently, another line of work has been devoted to strong convexification techniques in solving mixed-integer problems for power systems~\cite{Ramtin17, Fattahi2017, Fattahi17_2}.

%All of these methods are based on the linearization of disjunctive constraints in the OTS problem using the so-called big-M method.

%Although the direct current (DC) approximation of power systems is used as a widely acceptable approximation of alternative current (AC) power flow equations, \cite{Coffrin14} argues that in some scenarios, this approximation may be far from accurate in the context of OTS problem.. Recently, \cite{Kocuk15} proposed a mixed-integer second order conic relaxation of the OTS problem for AC systems.

In this work, the power flow equations are modeled using the well-known DC approximation, which is the backbone of the operation of power systems. Despite its shortcomings for the OTS in some cases \cite{Coffrin14}, the DC approximation is often considered very useful for increasing the reliability, performance, and market efficiency of power systems \cite{Hedman2011}. The OTS consists of disjunctive constraints that are bilinear and nonconvex in the original formulation. However, all of these constraints can be written in a linear form using the so-called big-$M$ or McCormick inequalities~\cite{Belotti11,McCormick76}. This formulation of OTS is referred to as the \textit{linearized OTS} in the sequel. A natural question arising in constructing the OTS formulation is: how can one find optimal values for the parameters of the big-$M$ or McCormick inequalities? An optimal choice for these parameters is important for two reasons: 1) they would result in stronger convex relaxations of the problem, and hence, fewer iterations in branch-and-bound or cutting-plane methods, and 2) a conservative choice of these parameters would cause numerical and convergence issues \cite{Wright98}.
Hedman {\it et al.} \cite{Hedman2009} point out that finding the optimal values for the parameters of the linearized OTS may be cumbersome, and, therefore, they impose restrictive constraints on the absolute angles of voltages at different buses at the expense of shrinking the feasible region. 

In this work, it is proven that finding the optimal values for the parameters of the MILP or MIQP formulations of the OTS using either big-$M$ or McCormick inequalities is NP-hard. Moreover, it is shown that there does not exist any polynomial-time algorithm to approximate these parameters within any constant factor, unless $P = NP$. This new result adds a new dimension to the difficulty of the OTS; not only is solving the OTS as a mixed-integer nonlinear program difficult, but finding a good linearized reformulation of this problem is NP-hard as well. 

%In practice, system operators consider flexibility only for a set of key transmission lines. 
In order to maintain the reliability and security of the system, often a set of transmission lines are considered as fixed and the flexibility in the network topology is limited to the remaining lines. An implicit requirement is that the network should always remain connected in order to prevent islanding. One way to circumvent the islanding issue in the optimal transmission switching problem is to include additional security constraints in order to keep the underlying network connected at every feasible solution~\cite{Khodaei10, Ostrowski14}. However, this new set of constraints would lead to the over-complication of an already difficult problem. Therefore, in practice, many energy corporations, such as PJM and Exelon, consider only a selected subset of transmission lines as flexible assets in their network~\cite{PJM2008, Exelon08}. In this paper, it is proven that the OTS with a connected spanning fixed subnetwork is still NP-hard but one can find non-conservative values for the parameters of the big-$M$ or McCormick inequalities in the linearized OTS. In particular, a simple bound strengthening method is presented to strengthen the linearized formulation of the OTS. 
%The proposed method is based on solving a simple shortest path problem on a weighted graph. 
This method can be integrated as a preprocessing step into any numerical solver for the OTS. Despite its simplicity, it is shown through extensive case studies on the IEEE 118-bus system and different Polish networks that the incorporation of the proposed bound strengthening method leads to substantial speedup in the runtime of the solver. 
%We perform our simulations on OTS problem with both linear and quadratic objective functions that are commonly used for minimizing the system loss and operational costs, respectively.

\section{Problem Formulation}

Consider a power network with $n_b$ buses, $n_g$ generators, and $n_l$ lines. This network can be represented by a graph, denoted by $G(\mathcal{B},\mathcal{L})$, where $\mathcal{B}$ is the set of buses indexed from 1 to $n_b$ and $\mathcal{L}$ is the set of lines indexed as $(i,j)$ to represent a connection between buses $i$ and $j$. In order to streamline the presentation, we assume an arbitrary direction for each line of the power system. Define $\mathcal N_l^+(i)$ as the set of endpoints of the outgoing lines at bus $i$. In particular, $\mathcal N_l^+(i)$ is defined as $\{j\in \mathcal{B}|(i,j)\in\mathcal{L}\}$. Similarly,   define $\mathcal N_l^-(i)$ as the set of endpoints of the incoming lines at bus $i$. In particular, we have  $\mathcal N_l^-(i) \triangleq \{j\in \mathcal{B}|(j,i)\in\mathcal{L}\}$. Denote $\mathcal{G} = \{1,2,...,n_g\}$ as the set of generators in the system. Furthermore, let $\mathcal{N}_g(i)$ be the indices of generators that are connected to bus $i$. Note that $\mathcal N_g(i)$ may be empty for a bus $i$. The variable $p_i$ corresponds to the active-power production of generator $i\in\mathcal{G}$ and the variable $\theta_i$ is the voltage angle at bus $i\in\mathcal{B}$. For every $(i,j)\in\mathcal{L}$, the variable $f_{ij}$ denotes the active flow from bus $i$ to bus $j$. Consider the set of lines $\mathcal{S}\subseteq\mathcal{L}$ that are equipped with on/off switches and define the decision variable $x_{ij}$ for every $(i,j)\in\mathcal{S}$ as the status of the line $(i,j)$. Let $n_s$ denote the cardinality of this set. We refer to the lines  belonging to $\mathcal{S}$ as \textit{flexible} lines and the remaining lines as \textit{fixed} lines. Notice that the decision variables $p_i$, $\theta_i$, and $f_{ij}$ are continuous, whereas $x_{ij}$ is binary. For simplicity of notation, define the  variable vectors
\begin{align}
	& \mathbf{p} \triangleq [p_1, p_2, ..., p_{n_g}]^\top,
	 &&\!\!\!\! \Theta \triangleq [\theta_1, \theta_2, ..., \theta_{n_b}]^\top,\nonumber\\
	& \mathbf{f} \triangleq [f_{i_1 j_1}, f_{i_2 j_2},...,f_{i_{n_l} j_{n_l}}]^\top,
	 &&\!\!\!\!\mathbf{x} \triangleq [x_{i_1 j_1}, x_{i_2 j_2},...,x_{i_{n_s} j_{n_s}}]^\top,
\end{align}
where the lines in $\mathcal{L}$ are labeled as $(i_1,j_1),...,  (i_{n_l},j_{n_l})$ such that the first $n_s$ lines denote the members of $\mathcal{S}$.
The objective function of the OTS is defined as $\sum_{i\in\mathcal{G}}g_i(p_i)$, where $g_i(p_i)$ takes the  quadratic form
\begin{align}\label{c_pi}
g_i(p_i) = a_i\times p_i^2+b_i\times p_i + c_i.
\end{align}
with $a_i\not= 0$ or the linear form
\begin{align}\label{c_pi2}
g_i(p_i) = b_i\times p_i + c_i.
\end{align}
for some numbers $a_i,b_i, c_i\geq 0$. In this paper, we consider both quadratic and linear objective functions, which may correspond to system loss and operational cost of generators. Every in-operation power system must satisfy operational constraints arising from  physical and security limitations. The physical limitations include the unit and line capacities. Furthermore, the power system must satisfy the power balance equations. On the security side, there may be a cardinality constraint on the maximum number of flexible lines that  can be switched off in order to avoid endangering the reliable operation of the system. Let the vector $\mathbf{d} = [d_1, d_2, ....,d_{n_b}]^\top$ collect the set of demands at all buses. Moreover, define $p_i^{\min}$ and $p_i^{\max}$ as the lower and upper bounds on the production level of generator $i$, and $f_{ij}^{\max}$ as the capacity of line $(i,j)\in\mathcal{L}$. Each line $(i,j)\in\mathcal{L}$ is associated with susceptance $B_{ij}$.

Using the above notations, the OTS is formulated as the following mixed-integer nonlinear problem:
\begin{subequations}\label{ots}
	\begin{align}
	\underset{
		\mathbf{f}, \mathbf{x}, \mathbf{\Theta}, \mathbf{p}
	}{\text{minimize}} \ \ \ \
	\!\!\!\!& 
	\underset{
		i\in\mathcal{G}
	}{\sum}
	g_i(p_i) \\
	\mathrm{s.t.} \quad\quad\quad\ \ \ \  x_{ij}&\in\{0,1\},&& \forall (i,j)\in\mathcal{S} \label{const1}\\
	p^{\min}_{k} &\leq p_{k}\leq p^{\max}_{k}, && \forall k\in\mathcal{G}\label{const2}\\
	-f^{\max}_{ij} x_{ij}&\leq f_{ij}\leq f^{\max}_{ij} x_{ij}, \label{const3} && \forall (i,j)\in\mathcal{S}\\
	-f^{\max}_{ij}&\leq f_{ij}\leq f^{\max}_{ij}, \label{const4}  && \forall (i,j)\in\mathcal{L}\backslash\mathcal{S}\\
	B_{ij}(\theta_i-\theta_j) x_{ij} &= f_{ij},&& \forall (i,j)\in\mathcal{S}\label{const5}\\
	B_{ij}(\theta_i-\theta_j) &= f_{ij},&&\forall (i,j)\in\mathcal{L}\backslash\mathcal{S}\label{const6}\\
	\!\!\sum_{k\in\mathcal{N}_g(i)}\!\! p_k\!-\!d_i &=\!\!\!\!\! \sum_{j\in\mathcal{N}^+_l(i)}\!\!f_{ij}-\!\!\!\!\! \sum_{j\in\mathcal{N}^-_l(i)}\!\!f_{ji},\label{const7}\!\!\!\!\!&&  \forall i\in\mathcal{B}\\
	\sum_{(i,j)\in\mathcal{S}}x_{ij}&\geq r\label{const9},
	\end{align}
\end{subequations}
	where
	\begin{itemize}
		\item[-] \eqref{const1} states that the status of each flexible line must be binary;
		\item[-] \eqref{const2} imposes lower and upper bounds on the production level of generating units;
		\item[-] \eqref{const3} and \eqref{const4} state that the flow over a flexible or fixed line must be within the line capacities when its switch is on, and it should be zero otherwise;  
		\item[-] \eqref{const5} and \eqref{const6} relate the flow over each line to the voltage angles of the two endpoints of the line if it is in service, and it sets the flow to zero otherwise;
		\item[-] \eqref{const7} requires that the power balance equation be satisfied at every bus;
		\item[-] \eqref{const9} states that at least $r$ flexible lines must be switched on.
	\end{itemize}
	
	Define $\mathcal{F}$ as the feasible region of~\eqref{ots}, i.e., the set of $\{\mathbf{f}, \mathbf{x}, \mathbf{\Theta}, \mathbf{p}\}$ satisfying~\eqref{const1}-~\eqref{const9}.
	
Due to space restrictions, we consider only one time slot of the system operation. However, the techniques developed in this paper can also be used for the OTS over multiple time slots with coupling constraints, such as ramping limits on the productions of the generators. As another generalization, one can consider a combined unit commitment and optimal transmission switching problem, as formulated in Hedman {\it et al.} \cite{Hedman2010}.
	%Although our developed methodology can be readily applied to these generalizations, our focus is mainly on one time slot OTS problem due to space restrictions.
	Henceforth, the term ``optimal solution" refers to a globally optimal solution rather than a locally optimal solution. 
	
	To streamline the presentation, the proofs of the lemmas and theorems are moved to~\cite{Salar17}.
	
	\section{Linearization of OTS}\label{sec:linear}
	The aforementioned formulation of the OTS belongs to the class of mixed-integer nonlinear programs. The nonlinearity of this optimization problem is, in part, caused by the multiplication of the binary variable $x_{ij}$ and the continuous variables $\theta_i$ and $\theta_j$ in \eqref{const5}. However, since this nonlinear constraint has a disjunctive nature, one can use the big-$M$ or McCormick reformulation technique to formulate it in a linear way. First, we consider the big-$M$ method, and then show that the same result holds for the McCormick reformulation scheme in the OTS. One can re-write \eqref{const5}  for each flexible line $(i,j)$ in the  form
	\begin{equation}\label{bigM}
		B_{ij}(\theta_i-\theta_j)-M_{ij}(1-x_{ij})\leq f_{ij}\leq B_{ij}(\theta_i-\theta_j)+M_{ij}(1-x_{ij})
	\end{equation} 
	for a \textit{large enough} scalar $M_{ij}$, which results in the {\it linearized OTS formulation}. The above inequality implies that if $x_{ij}$ equals 1, then the line is in service and needs to satisfy the physical constraint $f_{ij} = B_{ij}(\theta_i-\theta_j)$. On the other hand, if $x_{ij}$ equals 0, then \eqref{bigM} (and hence \eqref{const5}) is redundant as it is dominated by \eqref{const3}.
	The term ``large enough'' for $M_{ij}$ is ambiguous, and indeed the design of an effective $M_{ij}$ is a challenging task that will be studied below. 
	\begin{definition}
		For every $(i,j)\in\mathcal{S}$, it is said that $M_{ij}$ is feasible for the OTS if it preserves the equivalence between \eqref{bigM} and \eqref{const5} in the OTS. The smallest feasible $M_{ij}$ is denoted by $M_{ij}^{\mathrm{opt}}$.
	\end{definition}

\begin{remark}
	{\rm Note that the value of $M_{ij}^{\mathrm{opt}}$ is \textit{independent} of the values of $M_{rl},$ for $(r,l)\in\mathcal{S}\backslash (i,j),$ in the linearized OTS formulation, as long as they are chosen to be feasible. In other words, given an instance of the OTS, the value of $M_{ij}^{\mathrm{opt}}$ is the same if $M_{rl}$ satisfies $M_{rl}\geq M_{rl}^{\mathrm{opt}}$ for every $(r,l)\in\mathcal{S}\backslash (i,j)$.}
\end{remark}

The problem under investigation in this section is the following:
Given an instance of OTS, is there an efficient algorithm to compute $M_{ij}^{\mathrm{opt}}$ or a good approximation of that for every $(i,j)\in\mathcal{S}$?
%Note that one can pick conservative values for $M_{ij}$ to  make \eqref{bigM} equivalent to \eqref{const5}. For instance, it can be easily verified that, for every flexible line $(i,j)$, the number $\sum_{(r,k)\in\mathcal{L}}f^{\max}_{rk}/B_{rk}$ is a feasible choice for $M_{ij}$. 
It is desirable to find the smallest feasible values for every $M_{ij}, (i,j)\in\mathcal{S},$ in~\eqref{bigM}  because of two reasons:
\begin{itemize}
	\item[1.] Commonly used methods for solving MILP or MIQP problems, such as cutting-plane and branch-and-bound algorithms, are based on iterative convex relaxations of the constraints. Therefore, while a sufficiently large value for $M_{ij}$ does not change the feasible region of the OTS after replacing \eqref{const5} with \eqref{bigM}, it may have a significant impact on the feasible region of its convex relaxation.  Small values for $M_{ij}$ yield stronger convex relaxations with smaller feasible sets.
	\item[2.] Large values for $M_{ij}$ may cause numerical issues for convex relaxation solvers. 
\end{itemize}

For every $(i,j)\in\mathcal{S}$, define $\mathcal{F}_{ij}$ as the set of all points $\{\mathbf{f}, \mathbf{x}, \mathbf{\Theta}, \mathbf{p}\}\in\mathcal{F}$ such that $x_{ij} = 0$. 
\begin{lemma}\label{Mopt}
	The equation
	\begin{equation}
		M_{ij}^{\mathrm{opt}} = B_{ij}\times \max_{\{\mathbf{f}, \mathbf{x}, \mathbf{\Theta}, \mathbf{p}\}\in\mathcal{F}_{ij}}\{|\theta_i-\theta_j|\}
	\end{equation}
	holds for every flexible line $(i,j)\in\mathcal{S}$.
\end{lemma}
%\begin{proof}
%	Consider a number $M_{ij}$ such that $M_{ij}\geq B_{ij}\times \max_{\mathcal{F}_{ij}}\{|\theta_i-\theta_j|\}$. Every $\{\mathbf{f}, \mathbf{x}, \mathbf{\Theta}, \mathbf{p}\}\in\mathcal{F}$ satisfies~\eqref{bigM} with the chosen $M_{ij}$ and hence, $M_{ij}$ is feasible. Now, assume that $M_{ij}< B_{ij}\times \max_{\mathcal{F}_{ij}}\{|\theta_i-\theta_j|\}$. It is straightforward to observe that there exists $\{\mathbf{f}, \mathbf{x}, \mathbf{\Theta}, \mathbf{p}\}\in\mathcal{F}$ for which~\eqref{bigM} is violated and, hence, $M_{ij}$ is not feasible. This implies that $M_{ij}^{\mathrm{opt}} = B_{ij}\times \max_{\mathcal{F}_{ij}}\{|\theta_i-\theta_j|\}$.
%\end{proof}
Due to Lemma~\ref{Mopt}, the problem of finding $M_{ij}^{\mathrm{opt}}$ for every $(i,j)\in\mathcal{S}$ reduces to finding the $\max_{\mathcal{F}_{ij}}\{|\theta_i-\theta_j|\}$.
\begin{remark}
	{\rm Note that, for a given $(i,j)\in\mathcal{S}$, the term $\max_{\mathcal{F}_{ij}}\{|\theta_i-\theta_j|\}$ is finite if and only if the buses $i$ and $j$ are connected for every feasible point in $\mathcal{F}_{ij}$. This means that the linearization of the OTS is well-defined if and only if the power network remains connected at every feasible solution in $\mathcal{F}_{ij}$ for all $(i,j)\in\mathcal{S}$.} 
\end{remark}
The next example illustrates a scenario where the $\max_{\mathcal{F}_{ij}}\{|\theta_i-\theta_j|\}$ is not finite. 
\begin{figure}
	\centering
	%\vspace*{-0.75cm}
	\includegraphics[width=0.8\columnwidth]{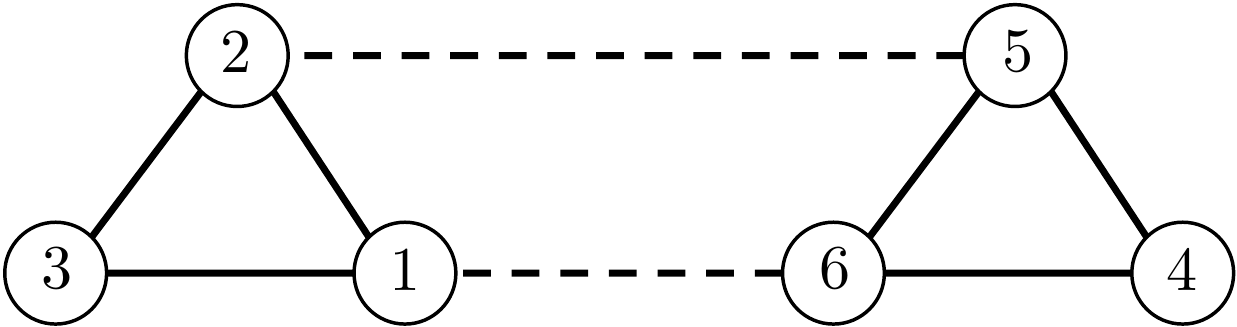}
	\caption{ \footnotesize This depicts the topology of the network in Example~\ref{exp1}. The solid and dashed edges denote the lines with ON and OFF switches, respectively.}
	\label{example}
\end{figure}

\begin{example}\label{exp1}
	{\rm Consider the network with 6 buses and 8 lines in Figure~\ref{example}. Assume that the network is decomposed into two disjoint components (known as islands) with the buses $\{1,2,3\}$ and $\{4,5,6\}$ at a feasible point $\{\mathbf{f}, \mathbf{x}, \mathbf{\Theta}, \mathbf{p}\}\in\mathcal{F}_{16}$. Define $\tilde{\mathbf{\Theta}}$ as $\tilde{\theta}_i = \theta_i$ for $i\in\{1,2,3\}$ and $\tilde{\theta}_i = \theta_i+\tau$ for $i\in\{4,5,6\},$ where $\tau$ is an arbitrary scalar. It can be verified that $\{\mathbf{f}, \mathbf{x}, \tilde{\mathbf{\Theta}}, \mathbf{p}\}\in\mathcal{F}_{16}$ for every $\tau$. Furthermore, $\tilde{\theta}_6-\tilde{\theta}_1 = {\theta}_6-{\theta}_1+\tau$, which implies that $\max_{\mathcal{F}_{16}}\{|\tilde{\theta}_6-\tilde{\theta}_1|\}\rightarrow +\infty$ as $\tau\rightarrow +\infty$.}
\end{example}

To avoid unbounded values for $M^{\mathrm{opt}}_{ij}$, the existence of a connected spanning subnetwork connecting all the nodes in the network with fixed lines will be assumed in the next section.
In what follows, it will be shown that, even if $\max_{\mathcal{F}_{ij}}\{|\theta_i-\theta_j|\}$ is bounded for every $(i,j)\in\mathcal{S}$, one cannot devise an algorithm that efficiently finds $\max_{\mathcal{F}_{ij}}\{|\theta_i-\theta_j|\}$ because it amounts to an NP-hard problem. Furthermore, the impossibility of any constant factor approximation of $\max_{\mathcal{F}_{ij}}\{|\theta_i-\theta_j|\}$ in the linearized OTS is proven.

%	\begin{definition}
%		Given an instance of OTS problem and a value $M^*_{ij}$ for the line $(i,j)\mathcal{S}$, the \textit{feasibility version} of OTS problem is defined as determining whether $M^*_{ij}\in\mathcal{M}_{ij}$.
%	\end{definition}
	\begin{theorem}\label{thm1}
		Consider an instance of the OTS and select a flexible line $(i,j)\in\mathcal{S}$. Unless $P=NP$, it holds that:
		\vspace{-0cm}
		\begin{itemize}
			\item[-] (NP-hardness) there is no polynomial-time algorithm for finding $\max_{\mathcal{F}_{ij}}\{|\theta_i-\theta_j|\}$;
			\item[-] (Inapproximability) there is no polynomial-time constant-factor approximation algorithm for finding $\max_{\mathcal{F}_{ij}}\{|\theta_i-\theta_j|\}$.
		\end{itemize}
	\end{theorem} 
	
%	\begin{proof}
%		The proof is provided in the Appendix.
%	\end{proof}
Theorem~\ref{thm1} together with Lemma~\ref{Mopt} implies that finding $M_{ij}^{\mathrm{opt}}$ is both NP-hard and inapproximable within any constant factor, hence providing a negative answer to the question raised in this section.
\begin{remark}
	{\rm The decision version of the OTS is known to be NP-complete \cite{Burak2016}. One may speculate that the NP-hardness of finding the best $M_{ij}$ for every $(i,j)\in\mathcal{S}$ may follow directly from that result. However, notice that there are some well-known problems with disjunctive constraints, such as the minimization of total tardiness on a single machine, which are known to be NP-hard \cite{Du1990} and yet there are efficient methods to find the optimal parameters of their big-$M$ reformulation \cite{Herr2014}. Theorem \ref{thm1} shows that not only is finding the best $M_{ij}$ for the OTS NP-hard, but one cannot hope for obtaining a strong linearized reformulation of the problem based on the big-$M$ method.}
\end{remark}

%\subsection{McCormick Inequalities}

Note that one may choose to use McCormick inequalities ~\cite{McCormick76} instead of the big-$M$ method to obtain a linear reformulation of the bilinear constraint~\eqref{const5}. In what follows, it will be shown that the complexity of finding the optimal parameters of McCormick inequalities is the same as those in the big-$M$ method for the OTS. The McCormick inequalities can be written in the following form for a flexible line $(i,j)$:
\begin{subequations}\label{MC}
	\begin{align}
	& f_{ij}\leq u_{ij|x_{ij} = 1}x_{ij},\\
	& f_{ij}\geq l_{ij|x_{ij} = 1}x_{ij},\\
	& f_{ij}\leq B_{ij}(\theta_i-\theta_j)-l_{ij|x_{ij} = 0}x_{ij},\\
	& f_{ij}\geq B_{ij}(\theta_i-\theta_j)-u_{ij|x_{ij} = 0}x_{ij},
	\end{align}
\end{subequations}
where $u_{ij|x_{ij} = 1}$ and $l_{ij|x_{ij} = 1}$ are the respective upper and lower bounds for $B_{ij}(\theta_i-\theta_j)$ in the case where the line $(i,j)$ is in service. Similarly, $u_{ij|x_{ij} = 0}$ and $l_{ij|x_{ij} = 0}$ are the respective upper and lower bounds for $B_{ij}(\theta_i-\theta_j)$ when the switch for the flexible line $(i,j)$ is off. It can be verified that the following equalities hold:
\begin{subequations}\label{MCbound}
	\begin{align}
	& u_{ij|x_{ij} = 1} = f^{\max}_{ij}, &&\!\!\!\!l_{ij|x_{ij} = 1} = -f^{\max}_{ij},\\
	& u_{ij|x_{ij} = 0} \!=\! B_{ij}\!\times\!\max_{\mathcal{F}_{ij}}\{\theta_i\!-\!\theta_j\}, && \!\!\!\!l_{ij|x_{ij} = 0} \!=\! B_{ij}\!\times\!\min_{\mathcal{F}_{ij}}\{\theta_i\!-\!\theta_j\}.
	\end{align}
\end{subequations}
Therefore, Theorem~\ref{thm1} immediately results in the NP-hardness and inapproximability of the pair $(l_{ij|x_{ij} = 0}, u_{ij|x_{ij} = 0})$.

%Although finding a good approximation of $M^{\mathrm{opt}}_{ij}$ is a difficult task in general, one can find a non-conservative value for $M_{ij}$ in polynomial time for power systems whose fixed lines form a connected spanning subgraph.

\section{OTS with a Fixed Connected Spanning Subgraph}\label{sec:fixed}

In this section, we consider a power system with the property that the set of fixed lines contains a connected spanning tree of the power system. The objective is to show that a non-trivial upper bound on $M_{ij}$ can be efficiently derived by solving a shortest path problem. Furthermore, it will be proven that this upper bound is tight in the sense that there exist instances of the OTS with a fixed connected spanning subgraph for which this upper bound equals $M^{\mathrm{opt}}_{ij}$. Before presenting this result, it is desirable to state that the OTS is hard to solve even under the assumption of a fixed connected spanning subgraph.
\begin{theorem}
	The OTS with a fixed connected spanning subgraph is NP-hard.
\end{theorem}
%\begin{proof}
%	The proof follows from a slight modification of the argument made in the proof of Theorem 3.1 in \cite{Burak2016}. 
%\end{proof}

Consider a feasible point $\{\mathbf{f}, \mathbf{x}, \mathbf{\Theta}, \mathbf{p}\}\in\mathcal{F}$. For any line $(i,j)\in\mathcal{L}$, we have 
\begin{equation}\label{eq12}
	B_{ij}(\theta_i-\theta_j) = B_{ij}\sum_{(r,l)\in \mathcal{P}_{ij}}(\theta_r-\theta_l),
\end{equation}
where $\mathcal{P}_{ij}$ is an arbitrary path from node $i$ to node $j$ in the fixed spanning connected subgraph of $G$. This implies that 
\begin{equation}\label{eq13}
	M^{\mathrm{opt}}_{ij} = B_{ij}\sum_{(r,l)\in \mathcal{P}_{ij}}(\theta^{\mathrm{opt}}_r-\theta^{\mathrm{opt}}_l)\leq B_{ij}\sum_{(r,l)\in\mathcal{P}_{ij}}\frac{f^{\max}_{rl}}{B_{rl}},
\end{equation}
where $\{\mathbf{f}^{\mathrm{opt}}, \Theta^{\mathrm{opt}},\mathbf{x}^{\mathrm{opt}}, \mathbf{p}^{\mathrm{opt}} \}\in\arg\max_{\mathcal{F}}\{|\theta_i-\theta_j|\}$. Note that~\eqref{eq13} holds for every path $\mathcal{P}_{ij}$ in the fixed connected spanning subgraph of the network. We will use this observation in Theorem~\ref{thm3} to derive strong upper bounds for $M^{\mathrm{opt}}_{ij}$. Denote the undirected weighted subgraph induced by the fixed lines in the power system as $G_{\mathcal{I}}(\mathcal{B}_{\mathcal{I}}, \mathcal{W}_{\mathcal{I}})$, where $\mathcal{B}_{\mathcal{I}} = \mathcal{B}$ and $\mathcal{W}_{\mathcal{I}}$ is the set of all tuples $(i,j,w_{ij})$ such that $(i,j)\in\mathcal{L}\backslash\mathcal{S}$ and $w_{ij}$ is the weight corresponding to $(i,j)$ defined as $f^{\max}_{ij}/B_{ij}$. Let $\mathcal{P}_{\mathcal{I};ij}$ and $p_{\mathcal{I};ij}$ be the set of edges in a shortest simple
path between nodes $i$ and $j$ in $G_{\mathcal{I}}$ and its length, respectively. 
%Define $\mathcal{P}_{\mathcal{I};ij}$ and $p_{\mathcal{I};ij}$ as sets of edges in the shortest simple path between nodes $i$ and $j$ in $G_{\mathcal{I}}$ and its length, respectively. For simplicity of notation, we assume that the shortest path is unique. The results in the sequel also hold when the shortest path between nodes $i$ and $j$ in $G_{\mathcal{I}}$ is not unique.
%Furthermore, we set the weight of each edge $(i,j)\in\mathcal{L}_{\mathcal{I}}$ equal to $f^{\max}_{ij}/B_{ij}$ and let $\mathcal{W}_{\mathcal{I}}$ denote the set of all these weights.
\begin{theorem}\label{thm3}
	For every flexible line $(i,j)\in\mathcal{S}$, the inequality
	\begin{equation}
		M_{ij}^{\mathrm{opt}}\leq B_{ij}\times p_{\mathcal{I},ij}
	\end{equation}
	holds. Moreover, there exists an instance of the OTS for which this inequality is tight.
\end{theorem} 
%\begin{proof}
%%	Suppose that $\{\mathbf{f}^{\mathrm{opt}}, \Theta^{\mathrm{opt}},\mathbf{x}^{\mathrm{opt}}, \mathbf{p}^{\mathrm{opt}} \}\in\arg\max_{\mathcal{F}}\{|\theta_i-\theta_j|\}$. This implies that $M^{\mathrm{opt}}_{ij} = B_{ij}(\theta^{\mathrm{opt}}_i-\theta^{\mathrm{opt}}_j)$, which is equal to 
%%	\begin{equation}\label{PI}
%%		B_{ij}\sum_{(r,l)\in \mathcal{P}_{\mathcal{I},ij}}(\theta^{\mathrm{opt}}_r-\theta^{\mathrm{opt}}_l)
%%	\end{equation}
%%	Notice that the existence of $\mathcal{P}_{\mathcal{I},ij}$ is due to the presence of a fixed connected spanning subgraph. The equation \eqref{PI} yields that
%Based on~\eqref{eq13}, we have 
%	\begin{equation}\label{shortest}
%		M_{ij}^{\mathrm{opt}} \leq B_{ij}\!\sum_{(r,l)\in\mathcal{P}_{\mathcal{I},ij}}\!\frac{f^{\max}_{rl}}{B_{rl}} = B_{ij}\!\!\sum_{(r,l)\in\mathcal{P}_{\mathcal{I},ij}}\!\!w_{rl} = B_{ij}\times p_{\mathcal{I},ij}.
%	\end{equation}
%	Furthermore, a similar feasible solution that is derived in the second part of the proof of Theorem \ref{thm1} (provided in the Appendix) can be used to show the tightness of this upper bound.
%\end{proof}
 
Theorem~\ref{thm3} proposes a bound strengthening scheme for every flexible line in the OTS problem that can be carried out as a simple preprocessing step before solving the OTS using any branch-and-bound method. The algorithm for the proposed bound strengthening method is described in Algorithm~\ref{alg1}. 
\begin{algorithm}[]
	\KwData{$G_{\mathcal{I}}(\mathcal{B}_{\mathcal{I}}, \mathcal{W}_{\mathcal{I}})$ and $B = \{B_{ij}| (i,j)\in\mathcal{S}\}$}
	\KwResult{$M_{ij}$ for every $(i,j)\in\mathcal{S}$}
	\For{$(i,j)\in\mathcal{S}$}{
		find $p_{\mathcal{I};ij}$ using Dijkstra's algorithm\;
		$M_{ij}\gets B_{ij}\times p_{\mathcal{I};ij}$\;
	}
	\caption{Bound strengthening method for linearized OTS with fixed connected spanning subgraph}\label{alg1}
\end{algorithm}

The worst-case complexity of performing this preprocessing step is ${O}(n_s n_b^2)$ since it is equivalent to performing $n_s$ rounds of Dijkstra's algorithm on the weighted graph $G_{\mathcal{I}}$ (it can also be reduced to ${O}(n_s(n_l-n_s+n_b\log n_b))$ if the algorithm is implemented using a Fibonacci heap)~\cite{Algo09}. This preprocessing step can be processed in an offline fashion before realizing the demand in the system. The impact of this preprocessing step on the runtime of the solver will be demonstrated on different cases in Section~\ref{sec:num}. 

As mentioned in the Introduction, the existence of a fixed connected spanning subgraph in power systems is a practical assumption since  power  operators should guarantee the reliability of the system by ensuring the connectivity of the power network. Therefore, due to Theorem \ref{thm3}, one can design relatively small values for $M_{ij}$'s in order to strengthen the convex relaxation of OTS. 
%In the remainder of the paper, we will focus on the OTS problem with a fixed connected subgraph. However, it is worthwhile to mention that the proposed convex model can be readily applied to the OTS problem without assuming the existence of such fixed connected subgraph, as long as a relatively good upper bound on $M^*_{ij}$ is known \textit{a priori} for every flexible line $(i,j)$.

%\section{Quadratic Cost Function}\label{sec:quad}

Consider the cost function for the OTS. In practice, for the production planning problems of power systems, it is often the case that a quadratic objective function is used in order to better imitate the actual cost of production, specially for thermal generators~\cite{Wood12}. However, the nonlinearity introduced by a quadratic cost function makes the OTS hard to solve in general. 
%Due to this inherent difficulty in solving this MIQP, most of the existing literature on the OTS assumes a linear cost function~\cite{Fisher2008,Burak2016} for the production of electricity in all generators. 
The main challenge of solving the MIQP is the fact that the optimal solution of its continuous relaxation often lies in the interior or on the boundary of its relaxed feasible region which may be infeasible for the original MIQP (as opposed to the extreme point solutions in MILP). More precisely, even obtaining the convex hull of the feasible region is not enough to guarantee the exactness of such continuous relaxations, since the optimal solution of the relaxed problem does not always correspond to an extreme point in the convex hull if the objective function is quadratic. 
%On the other hand, branch-and-bound methods heavily rely on iterative continuous relaxations of the original problem. 
This would introduce non-fractional solutions for the binary variables of the problem in most of the iterations of branch-and-bound methods which often leads to a high number of iterations. 
One way to partially remedy this problem is to reformulate the problem by introducing auxiliary variables such that a new linear function is minimized and the old quadratic objective function is moved to the constraints. This guarantees that the continuous relaxation of the reformulated problem will obtain an optimal solution that is an extreme point of the relaxed feasible region. This is a key reason behind the success of different conic relaxation and strengthening methods in MIQP~\cite{Selim09, Alper07}. 
%In other words, different valid inequalities and bound strengthening methods will be more effective since they can partially describe the convex hull of the reformulated problem with a linear objective function. 

Assume that the objective function is quadratic in the form of $\sum_{i=1}^{n_g}g_{i}(p_i)$, where $g_i(p_i)$ is defined as~\eqref{c_pi}. Upon defining a new set of variables $t_i$ for $i\in\mathcal{G}$, one can reformulate the objective function as $\sum_{i=1}^{n_g}\tilde{g}_i(p_i,t_i)$ where
\begin{align}\label{c_pi_tilde}
\tilde{g}_i(p_i,t_i) = a_i\times t_i+b_i\times p_i + c_i.
\end{align}
subject to the additional convex constraints
\begin{equation}\label{quad}
p_i^2\leq t_i,\qquad \forall i\in\mathcal{G}
\end{equation}
This problem is referred to as \textit{modified formulation} of OTS henceforth.
%One can easily verify that the aforementioned modifications lead to an exact reformulation of the original OTS. However, it will be shown through different case studies that our proposed bound strengthening method is notably more effective when applied to the OTS with the linear objective function~\eqref{c_pi_tilde} and the additional set of constraints~\eqref{quad}. 

\section{Numerical Results}\label{sec:num}
In this section, case studies on different test cases are conducted to evaluate the effectiveness of the proposed preprocessing method in solving the OTS. All of the test cases are chosen from the publicly available MATPOWER package~\cite{Zimmerman11, Coffrin141}. The simulations are run on a laptop computer with an Intel Core i7 quad-core 2.50 GHz CPU and 16GB RAM. The results reported in this section are for a serial implementation in MATLAB using the CVX framework and the GUROBI 6.00 solver with the default settings. The relative optimality gap threshold is defined as
\begin{equation}\nonumber
\frac{z_{UB} - z_{LB}}{z_{UB}}\times 100,
\end{equation}
where $z_{UB}$ and $z_{LB}$ are the objective value corresponding to the best found feasible solution and the best found lower bound, respectively.  If the solver obtains a feasible solution for the OTS with the relative optimality gap of at most $0.1\%$ within a time limit (to be defined later), it is said that an optimal solution is found.

\subsection{Data Generation}
First, we study the IEEE 118-bus system. There are 185 lines in this test case. In all of the considered instances, a randomly generated connected spanning subgraph of the network with 120 fixed lines is chosen and the remaining lines are considered flexible. We compare the proposed preprocessing method with the case where $M_{ij}$ is chosen as $\sum_{(i,j)\in\mathcal{L}}f_{ij}^{\max}/B_{ij}$ for every $(i,j)\in\mathcal{S}$. This conservative value does not exploit the underlying structure of the network. To generate multiple instances of the OTS, the loads are multiplied by a load factor $\alpha$ chosen from the set $\{\alpha_1, \alpha_2,..., \alpha_k\}$. Furthermore, a uniform line rating is considered for all lines in the system. We examine both linear and quadratic cost functions and perform the following comparisons:

\begin{itemize}
	\item For the instances with a linear cost function, the total runtime is computed for the cases with and without the proposed bound strengthening method (denoted by \texttt{L-T} and \texttt{L-C}, respectively) for different load factors and cardinality lower bounds. 
	\item For the instances with a quadratic cost function, the runtime is computed for four different formulations: 1) the modified formulation with the proposed bound strengthening method (denoted by \texttt{Q-ET}), 2) the modified formulation without the proposed bound strengthening method (denoted by \texttt{Q-EC}), 3) the original formulation with the proposed bound strengthening method (denoted by \texttt{Q-OT}), and 4) the original formulation without the proposed bound strengthening method (denoted by \texttt{Q-OC}).
\end{itemize}

We also study six different large-scale Polish networks that are equipped with hundreds of switches. For each test case, a single load factor is considered for the OTS with linear and quadratic cost functions and the effect of the proposed bound strengthening method on the runtime is investigated. Similar to the IEEE 118-bus case, we fix a randomly chosen connected spanning subgraph of the network with fixed lines is considered. 

\subsection{IEEE 118-bus System}

\begin{figure*}
	\centering
	\subfloat[Running time vs. load factor]{\label{118_linear_lf}
		\includegraphics[width=0.97\columnwidth]{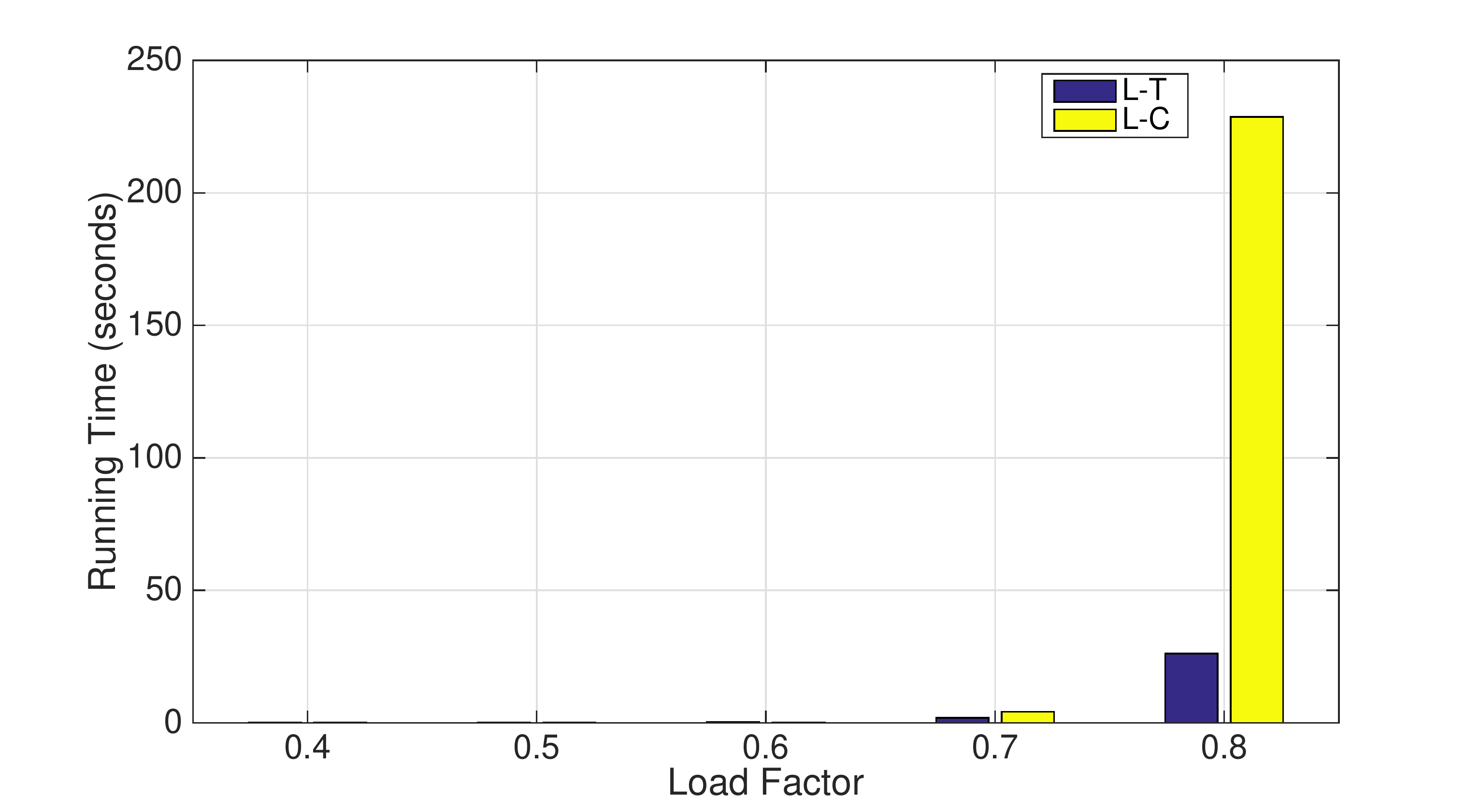}}
	\hspace{0.1cm}
	\subfloat[Running time vs. cardinality lower bound]{\label{118_linear_card}
		\includegraphics[width=0.97\columnwidth]{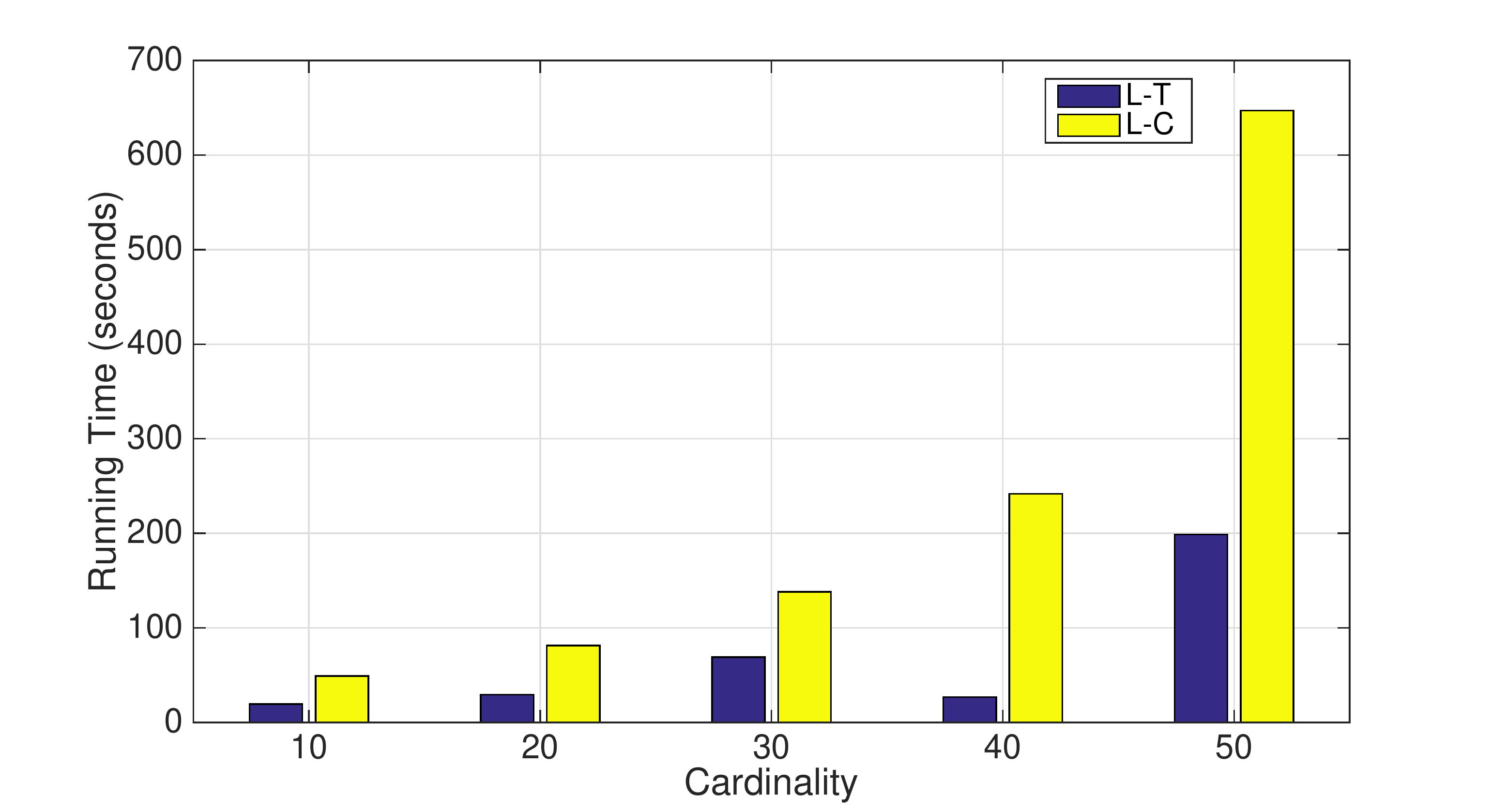}}
	\caption{ \footnotesize These plot show the runtime of different formulations of OTS with a linear cost function with respect to different load factors and cardinality lower bounds. \texttt{L-T} and \texttt{L-C} correspond to the original formulation with and without the bound strengthening method, respectively.}
	\vspace{-0.3cm}
\end{figure*}

In this subsection, the OTS problem is studied for the IEEE 118-bus system with 65 switches. Two types of cost functions are considered for this system:

\vspace{2mm}
\noindent{\bf Linear cost function:} Figure~\ref{118_linear_lf} shows the runtime with respect to the various load factors. For all of these experiments, the lower bound on the cardinality of the ON switches is set to 45, i.e. $r = 45$ in~\eqref{const9}. It can be observed that, for small values of the load factor, the OTS is relatively easy to solve with a linear cost function and the solver can easily find the optimal solution within a fraction of second with or without the bound strengthening method. On the other hand, as the load factor increases, the OTS becomes harder to solve and the proposed bound strengthening method has a significant impact on the runtime. In particular, when the load factor equals 0.8, the strengthened formulation of the OTS is solved 8.73 faster. 

In the second experiment, the performance of the solver is evaluated as a function of the lower bound on the number of the ON switches. As pointed out in~\cite{Burak2016}, the OTS becomes computationally hard to solve with a relatively large lower bound. This behavior is observed in Figure~\ref{118_linear_card}. However, note that the negative effect of increased lower bound diminishes when the bound strengthening step is performed. Specifically, the strengthened formulation is solved 2.66 times faster on average for the first two cardinality lower bounds (10 and 20) and 6.53 times faster on average for the last two cardinality lower bounds (40 and 50). 

\vspace{0.2cm}

\noindent{\bf Quadratic cost function:} When the cost function is quadratic, the runtime of the solver is drastically increased. Nevertheless, the modified formulation of the OTS combined with the proposed bound strengthening method reduces the runtime significantly.
\begin{figure*}
	\centering
	\subfloat[Running time vs. load factor]{\label{118_quad_lf}
		\includegraphics[width=0.97\columnwidth]{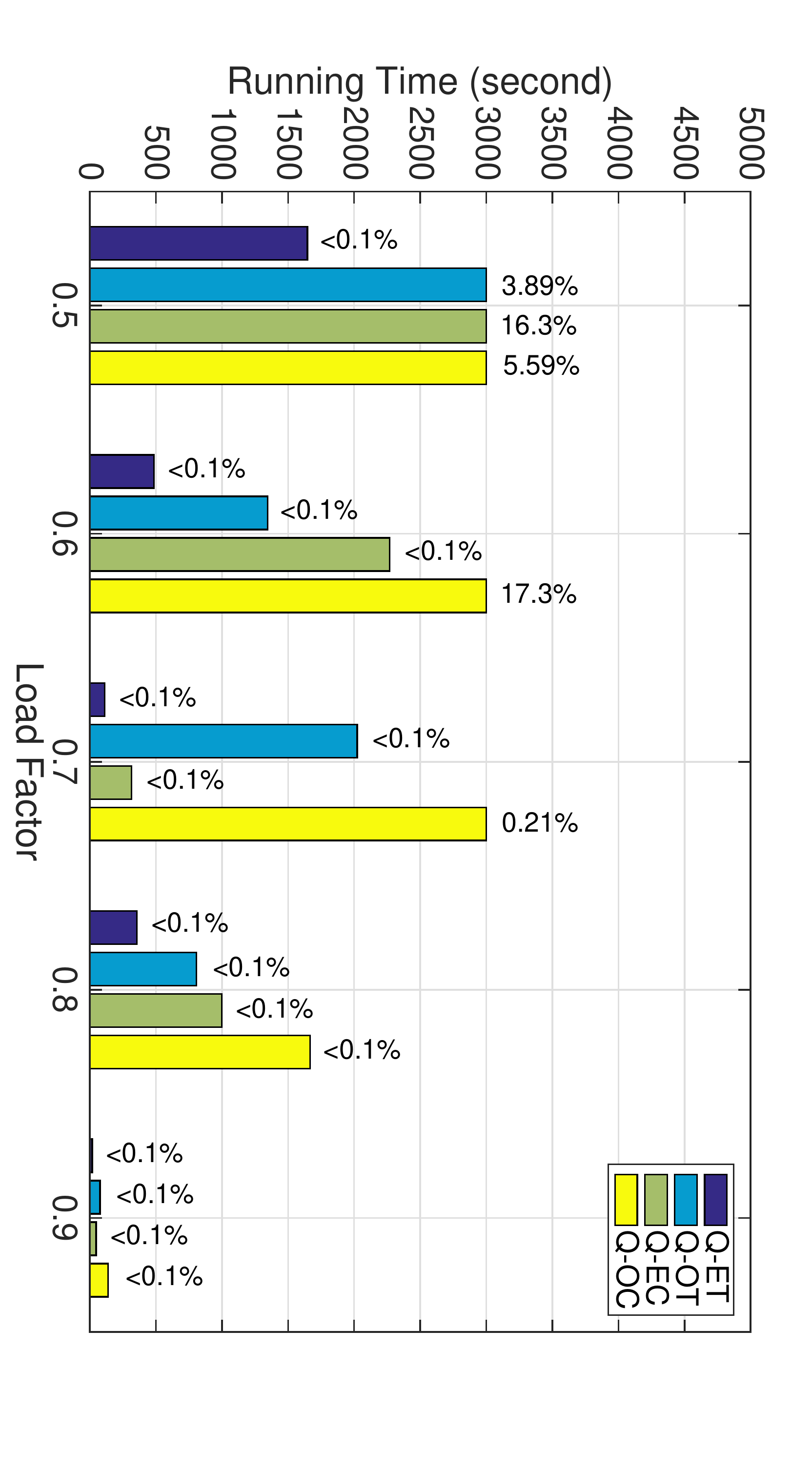}}
	\hspace{0.1cm}
	\subfloat[Running time vs. cardinality lower bound]{\label{118_quad_card}
		\includegraphics[width=0.97\columnwidth]{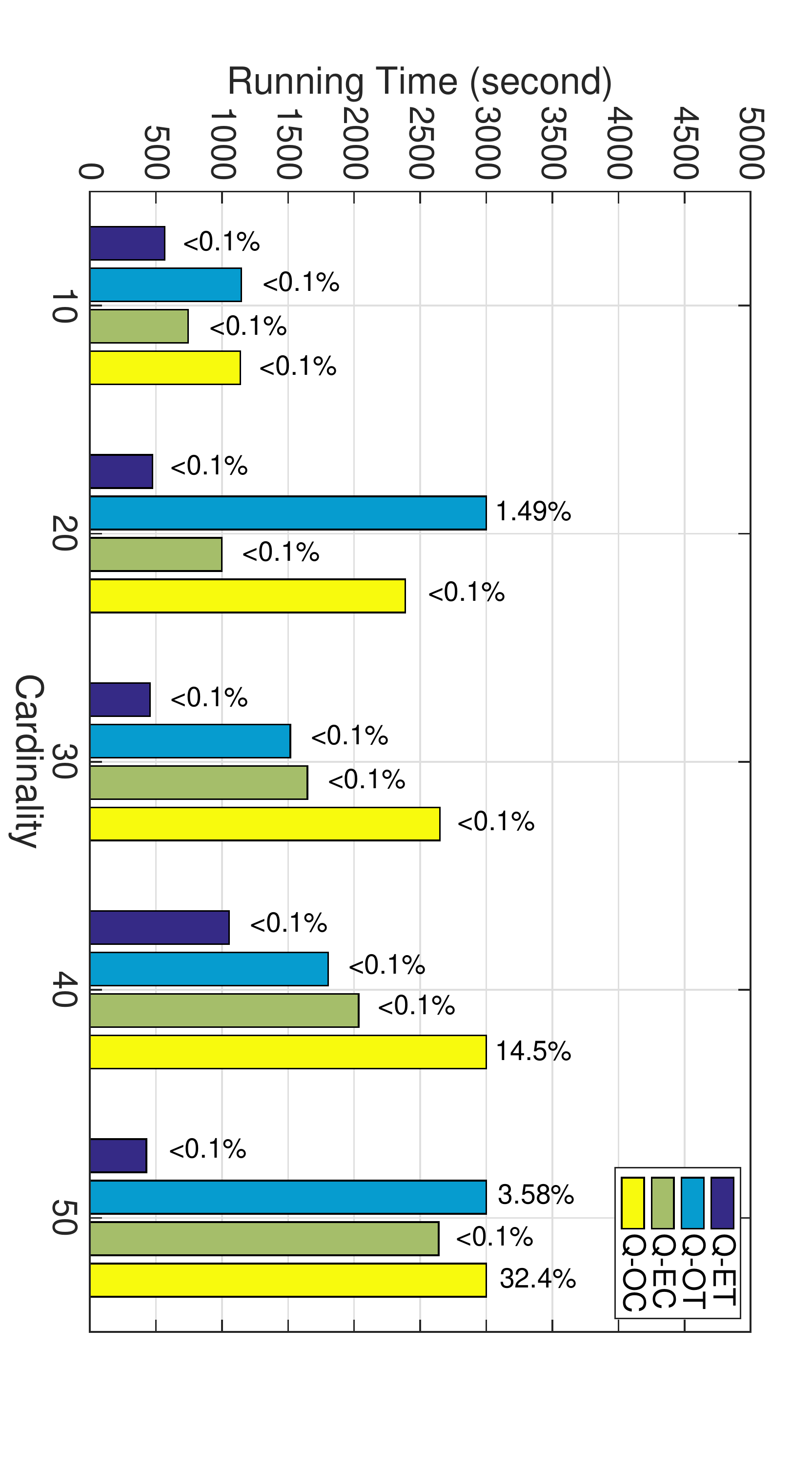}}
	\caption{ \footnotesize These plot show the runtime of different formulations of OTS with a quadratic cost function with respect to different load factors and cardinality lower bounds. \texttt{Q-ET}, \texttt{Q-EC}, \texttt{Q-OT}, and \texttt{Q-OC} correspond to the modified formulation with the proposed bound strengthening method, the modified formulation without the proposed bound strengthening method, the original formulation with the proposed bound strengthening method, and the original formulation without the proposed bound strengthening method, respectively.}
	\vspace{-0.3cm}
\end{figure*}
%\begin{figure}
%	\centering
%	%\vspace*{-0.75cm}
%	\includegraphics[width=1.05\columnwidth]{fig_lf.eps}
%	\caption{ \footnotesize The running time for different different formulation of OTS problem for IEEE 118-bus system with and without bound tightening method.}
%	\label{118_quad_card}
%\end{figure}
For all experiments, a time limit of $3,000$ seconds is imposed. For those instances that are not solved within the time limit, the relative optimality gap that is achieved by the solver at termination is reported. The runtime for different formulations of the OTS with respect to various load factors is depicted in Figure~\ref{118_quad_lf}. Similar to the previous case, the lower bound on the cardinality of the switches is set to 45 for different load factors. It can be observed that when the load factor equals 0.5, the solver can find the optimal solution within the time limit only for \texttt{Q-ET}. As the load factor increases, the average runtime decreases for all formulations. As it is clear from Figure~\ref{118_quad_lf}, \texttt{Q-ET} significantly outperforms other formulations for all load factors. Specifically, the runtime for \texttt{Q-ET} is at least 5.95, 2.96, and 13.58 times faster than \texttt{Q-OT}, \texttt{Q-EC}, and \texttt{Q-OC} on average, respectively. Notice that these values are the under-estimators of the actual speedups since the solver was terminated before finding the optimal solution in many cases. 

Next, consider the runtime for different formulations with respect to the change in the cardinality lower bound of ON switches. It can be observed that the solution times for \texttt{Q-OT}, \texttt{Q-EC}, and \texttt{Q-OC} increase as the lower bound increases. This observation supports the argument made in~\cite{Burak2016} suggesting that a large lower bound on the cardinality of the ON switches would make the OTS harder to solve in general. However, notice that the cardinality constraint has a minor effect on the runtime of \texttt{Q-ET}. 
Notice that \texttt{Q-OC} has the worst runtime on average among different settings of the load factor and cardinality lower bound. This implies that the proposed reformulation of the objective function together with the bound strengthening step is crucial to efficiently solve the OTS with a quadratic objective function.

\subsection{Polish Networks}

In this part, the proposed bound strengthening method is applied to solve the OTS for Polish networks. As for the 118-bus system, the runtime is evaluated for both linear and quadratic cost functions. In all of the simulations, the cardinality lower bound on the number of ON switches is set to 0. The number of flexible lines varies from 70 to 400. The time limit is chosen as $14,400$ seconds (4 hours) for the solver. If the time limit is reached, the optimality gap of the best found feasible solution (if one exists) is reported.
For the test cases with a quadratic cost function, only the modified formulation of the problem is considered because it significantly outperforms the original formulation.

\begin{table*}
	\caption{ \footnotesize This table shows the runtime of the solver with and without bound strengthening for Polish networks with linear and quadratic cost functions. The superscript $^*$ corresponds to the cases where the solver is terminated before finding the optimal solution due to the time limit.}
	\begin{center}
		\begin{tabular}{ c||c|c|c||c|c||c|c||c }
			& & & & \multicolumn{2}{|c||}{With Bound Tightening} & \multicolumn{2}{|c||}{Without Bound Tightening}&  \\
			\hline
			Cases & Cost Function & $\#$ Continuous & $\#$ Binary & Time & Optgap & Time & Optgap & Speedup\\
			\hline
			\multirow{2}{4em}{3120sp} & Linear & $3466$ & $70$ & $477$ & $<0.1\%$ & $3,623$ & $<0.1\%$ & $7.60$\\
			& Quadratic & $3466$ & $70$ & $2,900$ & $<0.1\%$ & $14,400$ & $0.12\%$ & $4.97^*$\\
			\hline
			\multirow{2}{4em}{2383wp} & Linear & $2789$ &$80$ & $418$ & $<0.1\%$ & $931$ & $<0.1\%$ & $2.23$ \\
			& Quadratic & $2789$ & $80$ & $252$ & $<0.1\%$ & $3,960$ & $<0.1\%$ & $15.71$\\
			\hline
			\multirow{2}{4em}{2736sp} & Linear & $3105$ &$100$ & $1,942$ & $<0.1\%$ & $2,508$ & $<0.1\%$ & $1.29$ \\
			& Quadratic & $3105$ & $100$ & $2,060$ & $<0.1\%$ & $3,417$ & $<0.1\%$ & $1.66$\\
			\hline
			\multirow{2}{4em}{3012wp} & Linear & $3516$ &$120$ & $2,447$ & $<0.1\%$ & $14,400$ & $0.11\%$ & $5.88^*$ \\
			& Quadratic &$3516$ &$120$ & $2,570$ & $<0.1\%$ & $14,400$ & $0.11\%$ & $5.60^*$\\
			\hline
			\multirow{2}{4em}{3375wp} & Linear & $4053$&$200$ & $98$ & $<0.1\%$ & $77$ & $<0.1\%$ & $0.79$\\
			& Quadratic & $4053$&$200$ & $4,301$ & $<0.1\%$ & $14,400$& -- & --\\
			\hline
			\multirow{2}{4em}{2746wop} & Linear & $3576$&$400$ & $17$ & $<0.1\%$ & $118$ & $<0.1\%$ & $6.94$\\
			& Quadratic &$3576$ &$400$ & $182$ & $<0.1\%$ & $3,523$ & $<0.1\%$ & $19.36$\\
			\hline
			\multirow{1}{4em}{\bf Average} & & & & {$\mathbf{1,472}$} & & $\mathbf{6,313}$ & & $\mathbf{6.55^*}$ \\
		\end{tabular}
	\end{center}\label{Polish}
\end{table*}

Table~\ref{Polish} reports the computational improvements when the bound strengthening method is incorporated into the formulation as a preprocessing step. This table includes the following columns:
\begin{itemize}
	\item $\#$ continuous: The number of continuous variables in the system;
	\item $\#$ binary: The number of binary variables corresponding to the flexible lines in the system;
	\item Time: The runtime (in seconds) for solving the OTS with and without bound strengthening method within the time limit;
	\item Optgap: The relative optimality gap within the time limit. Note that the solver is terminated when this gap is less than $0.1\%$;
	\item Speedup: The speedup in the runtime when the proposed bound strengthening method is used as a preprocessing step.
\end{itemize}

It can be observed from Table~\ref{Polish} that the presented bound strengthening method can notably reduce the computation time. In particular, the solver can be up to 19.36 times faster if the bound strengthening method is used to strengthen the formulation. Moreover, on average (excluding the case 3375wp with a quadratic cost function), the solution time is at least 6.55 times faster if the bound strengthening method is performed prior to solving the problem. For the case 3375wp with a quadratic cost function, the solver cannot obtain a feasible solution in $14,400$ seconds without bound strengthening. However, the solver can find an optimal solution within $4,301$ seconds after performing the proposed preprocessing step.

\section{Conclusion}

Finding an optimal topology of a power system subject to  operational and security constraints is a daunting task. In this problem, certain lines are fixed/uncontrollable, whereas the remaining ones could be controlled via on/off switches. The objective is to co-optimize the topology of the grid and the parameters of the system (e.g., generator outputs). Common techniques for solving this problem are mostly based on mixed-integer linear or quadratic reformulations using the big-$M$ or McCormick inequalities followed by iterative methods, such as branch-and-bound or cutting-plane algorithms. The performance of these methods partly relies on the strength of the convex relaxation of these reformulations. In this paper, it is shown that finding the optimal parameters of a linear or convex reformulation based on big-$M$ or McCormick inequalities is NP-hard. Furthermore, the inapproximability of these parameters up to any constant factor is proven. Despite the negative results on the complexity of the problem, a simple bound strengthening method is developed to significantly strengthen mixed-integer reformulations of the OTS, provided that there exists a connected spanning subgraph of the network with fixed lines. This bound strengthening method can be used as a preprocessing step even in an offline fashion, before forecasting the demand in the system. Through extensive computational experiments, it is verified that this simple preprocessing technique can significantly improve the runtime of the mixed-integer solvers in different test cases, including the IEEE 118-bus system and Polish networks.

\vspace{2mm}
\bibliographystyle{IEEETran}
\bibliography{reference}

\end{document}